\documentclass[12pt,oneside,english]{amsart}
\usepackage[T1]{fontenc}
\usepackage[latin1]{inputenc}
\usepackage{geometry}
\usepackage{setspace}
\usepackage{amssymb}

\makeatletter
 \theoremstyle{plain}
\newtheorem{thm}{Theorem}[section]
  \theoremstyle{plain}
  \newtheorem{lem}[thm]{Lemma}
  \theoremstyle{plain}
  \newtheorem{prop}[thm]{Proposition}

\numberwithin{equation}{section}

\usepackage{babel}
\makeatother
\begin{document}

\title{the asymptotic behavior of free additive convolution }

\author{hari bercovici and jiun-chau wang}

\begin{abstract}
We provide a new proof of the limit theorem for sums of free random
variables in a general infinitesimal triangular array. This result
was proved by Chistyakov and G\"{o}tze using subordination functions.
Our proof does not depend on subordination, and is close to the approach
used in the case of arrays with identically distributed rows \cite{BP in stable law}. 
\end{abstract}

\thanks{The first author was supported in part by a grant from the National
Science Foundation.}

\maketitle

\section{Introduction}

Given two probability measures $\mu$, $\nu$ on the real line $\mathbb{R}$,
we will denote by $\mu*\nu$ their classical convolution, and by $\mu\boxplus\nu$
their free additive convolution. Thus, $\mu*\nu$ is the distribution
of the sum $X+Y$, where $X$ and $Y$ are classically independent
random variables with distributions $\mu$ and $\nu$, respectively.
Analogously, $\mu\boxplus\nu$ is the distribution of $X+Y$, where
$X$ and $Y$ are freely independent random variables with distributions
$\mu$ and $\nu$. 

A triangular array $\{\mu_{nk}:n\geq1,\,1\leq k\leq k_{n}\}$ of probability
measures on $\mathbb{R}$ is said to be \emph{infinitesimal} if\[
\lim_{n\rightarrow\infty}\max_{1\leq k\leq k_{n}}\mu_{nk}(\{ t\in\mathbb{R}:\,\left|t\right|\geq\varepsilon\})=0,\]
for every $\varepsilon>0$. The classical limit distribution theory
for sums of independent random variables is concerned with the study
of the asymptotic behavior of the measures \[
\mu_{n}=\mu_{n1}*\mu_{n2}*\cdots*\mu_{nk_{n}}*\delta_{c_{n}},\qquad n\ge1,\]
where $\delta_{c_{n}}$ is the point mass at $c_{n}\in\mathbb{R}$.
Hin\v{c}in \cite{Hincin} proved that any weak limit of such a sequence
$\{\mu_{n}\}_{n=1}^{\infty}$ is an infinitely divisible measure.
Later Gnedenko (see \cite{GK} and \cite{Petrov}) found necessary
and sufficient conditions for the convergence of the sequence $\{\mu_{n}\}_{n=1}^{\infty}$
to a given infinitely divisible measure.

The analogous free convolutions\[
\nu_{n}=\mu_{n1}\boxplus\mu_{n2}\boxplus\cdots\boxplus\mu_{nk_{n}}\boxplus\delta_{c_{n}}\]
have been also the subject of several investigations. The first result
in this direction was an analogue of the central limit theorem proved
by Voiculescu \cite{V in symm}. Later, Pata \cite{Pata in CLT} proved
that the free central limit theorem holds precisely under the same
conditions as the classical central limit theorem. The analogue of
Hin\v{c}in's theorem, i.e. the fact that the possible weak limits
of $\nu_{n}$ are $\boxplus$-infinitely divisible, was proved in
\cite{BP in free hincin}. Then it was shown in \cite{BP in stable law}
that, in case $c_{n}=0$ and $\mu_{n1}=\mu_{n2}=\cdots=\mu_{nk_{n}}$,
the measures $\mu_{n}$ have a weak limit if and only if the measures
$\nu_{n}$ do. The correspondence between the limit of $\mu_{n}$
and the limit of $\nu_{n}$ was thoroughly studied in \cite{BT1,BT2}.

The result of \cite{BP in stable law} was extended in \cite{CG2}
to arbitrary arrays and centering constants $c_{n}$. The argument
in \cite{CG2} depends on two ingredients. The first is the fact that
the classical centering of the measures in an infinitesimal array
balances the real and the imaginary parts of the Cauchy transforms
of the measures. The second is the existence of subordination functions
as in \cite{V in free entropy,Biane,BB in semigroup}. 

We will provide a proof of the main result of \cite{CG2} which makes
no use of subordination functions, and is close to the argument of
\cite{BP in stable law}. This approach also works for multiplicative
free convolutions, as shown in \cite{BW in muli limit thm}.

The remainder of this paper is organized as follows. In Section 2
we describe the calculation of free convolution via Cauchy transform,
and we provide some useful approximation results. The proof of the
main result is in Section 3.

\section{Preliminaries}

Let $\mathcal{M}$ be the collection of all Borel probability measures
on $\mathbb{R}$. The free analogue of the Fourier transform was discovered
by Voiculescu \cite{V in addition} (see also \cite{Maassen} and
\cite{Bercovici and Voiculescu}). The details are as follows. Denote
by $\mathbb{C}^{+}=\{ z\in\mathbb{C}:\,\Im z>0\}$ the upper half-plane,
and set $\mathbb{C}^{-}=-\mathbb{\mathbb{C}}^{+}$. For a measure
$\mu\in\mathcal{M}$, one defines its \textit{Cauchy transform} $G_{\mu}:\:\mathbb{C}^{+}\rightarrow\mathbb{C}^{-}$
by 

\[
G_{\mu}(z)=\int_{-\infty}^{\infty}\frac{1}{z-t}\, d\mu(t),\qquad z\in\mathbb{C}^{+}.\]
Define the analytic function $F_{\mu}:\:\mathbb{C}^{+}\rightarrow\mathbb{C}^{+}$
by $F_{\mu}(z)=1/G_{\mu}(z)$. Given $\alpha,\beta>0$, set $\Gamma_{\alpha}=\{ z=x+iy\in\mathbb{C}^{+}:\,\left|x\right|<\alpha y\}$
and $\Gamma_{\alpha,\beta}=\{ z=x+iy\in\Gamma_{\alpha}:\, y>\beta\}$.
In \cite{Bercovici and Voiculescu} it is shown that $F_{\mu}(z)/z$
tends to 1 as $z\rightarrow\infty$ \textit{nontangentially} to $\mathbb{R}$
(i.e\textit{.,} $\left|z\right|\rightarrow\infty$ but $\Re z/\Im z$
stays bounded; in other words, $z\in\Gamma_{\alpha}$ for some $\alpha>0$),
and this implies that for every $\alpha>0$, there exists a $\beta=\beta(\alpha,\mu)>0$
such that $F_{\mu}$ has a left inverse $F_{\mu}^{-1}$ defined on
$\Gamma_{\alpha,\beta}$. The \textit{Voiculescu transform} of the
measure $\mu$ is then defined as\[
\phi_{\mu}(z)=F_{\mu}^{-1}(z)-z,\qquad z\in\Gamma_{\alpha,\beta}.\]
 We have $\Im\phi_{\mu}(z)\leq0$ for $z\in\Gamma_{\alpha,\beta}$,
and $\phi_{\mu}(z)=o(\left|z\right|)$ as $z\rightarrow\infty$ nontangentially. 

The most important property of the Voiculescu transform is that it
linearizes the free convolution. More precisely, if $\mu,\nu\in\mathcal{M}$
then\[
\phi_{\mu\boxplus\nu}(z)=\phi_{\mu}(z)+\phi_{\nu}(z),\]
for all $z$ in any truncated cone $\Gamma_{\alpha,\beta}$ where
all three functions involved are defined.

It was first noted in \cite{BP in free hincin} that for any given
cone $\Gamma_{\alpha,\beta}$, $\phi_{\mu}$ is defined on $\Gamma_{\alpha,\beta}$
as long as the measure $\mu$ puts most of its mass around the origin.

\begin{lem}
For every $\alpha,\beta>0$, there exists $\varepsilon>0$ with the
property that if $\mu\in\mathcal{M}$ is such that $\mu(\{ t\in\mathbb{R}:\,\left|t\right|\geq\varepsilon\})<\varepsilon$,
then $\phi_{\mu}$ is defined in $\Gamma_{\alpha,\beta}$.
\end{lem}
The following theorem from \cite{BP in LLN} translates weak convergence
of probability measures in terms of convergence properties of their
Voiculescu transforms. 

\begin{thm}
Given a sequence $\{\mu_{n}\}_{n=1}^{\infty}\subset\mathcal{M}$.
The following statements are equivalent:
\begin{enumerate}
\item $\mu_{n}$ converges weakly to a measure $\mu\in\mathcal{M}$ as $n\rightarrow\infty$; 
\item there exists a truncated cone $\Gamma_{\alpha,\beta}$ such that functions
$\phi_{\mu_{n}}$ and $\phi_{\mu}$ are defined in $\Gamma_{\alpha,\beta}$,
$\lim_{n\rightarrow\infty}\phi_{\mu_{n}}(iy)=\phi_{\mu}(iy)$ for
all $y>\beta$, and $\phi_{\mu_{n}}(iy)=o(y)$ uniformly in $n$ as
$y\rightarrow\infty$.
\end{enumerate}
\end{thm}
The following result is a reformulation of Propositions 2.6 and 2.7
in \cite{BP in stable law} which is more appropriate for our purposes. 

\begin{prop}
Given an infinitesimal array $\{\mu_{nk}:n\geq1,\,1\leq k\leq k_{n}\}\subset\mathcal{M}$,
and $\alpha,\beta>0$, the functions $\phi_{\mu_{nk}}$ are defined
in $\Gamma_{\alpha,\beta}$ for sufficiently large $n$, and \[
\phi_{\mu_{nk}}(z)=z^{2}\left[G_{\mu_{nk}}(z)-\frac{1}{z}\right](1+v_{nk}(z)),\]
 where the sequence \[
v_{n}(z)=\max_{1\leq k\leq k_{n}}\left|v_{nk}(z)\right|\]
satisfies $\lim_{n\rightarrow\infty}v_{n}(z)=0$ for all $z\in\Gamma_{\alpha,\beta}$,
and $v_{n}(z)=o(1)$ uniformly in $n$ as $z\rightarrow\infty$, $z\in\Gamma_{\alpha,\beta}$. 
\end{prop}
Recall that a measure $\nu\in\mathcal{M}$ is said to be \emph{$\boxplus$-infinitely
divisible} if for each $n\in\mathbb{N}$, there exists a measure $\mu_{n}\in\mathcal{M}$
such that\[
\nu=\underbrace{\mu_{n}\boxplus\mu_{n}\boxplus\cdots\boxplus\mu_{n}}_{n\,\text{times}}.\]
The notion of \emph{$*$-infinite divisibility} is defined similarly.
The well-known L\'{e}vy-Hin\v{c}in formula characterizes the $*$-infinitely
divisible measures in terms of their Fourier transform as follows:
a measure $\nu$ is $*$-infinitely divisible if and only if there
exist $\gamma\in\mathbb{R}$ and a finite positive Borel measure $\sigma$
on $\mathbb{R}$ such that the Fourier transform $\hat{\nu}$ is given
by \[
\hat{\nu}(t)=\exp\left[i\gamma t+\int_{-\infty}^{\infty}\left(e^{itx}-1-\frac{itx}{1+x^{2}}\right)\frac{1+x^{2}}{x^{2}}\, d\sigma(x)\right],\qquad t\in\mathbb{R}.\]
Here $\left(e^{itx}-1-\frac{itx}{1+x^{2}}\right)\frac{1+x^{2}}{x^{2}}$
is interpreted as $-t^{2}/2$ for $x=0$. We will denote by $\nu_{*}^{\gamma,\sigma}$
the $*$-infinitely divisible measure determined by $\gamma$ and
$\sigma$. The free analogue of the L\'{e}vy-Hin\v{c}in formula is
proved in \cite{Bercovici and Voiculescu}. A measure $\nu\in\mathcal{M}$
is $\boxplus$-infinitely divisible if and only if there exist $\gamma\in\mathbb{R}$
and a finite positive Borel measure $\sigma$ on $\mathbb{R}$ such
that\[
\phi_{\nu}(z)=\gamma+\int_{-\infty}^{\infty}\frac{1+tz}{z-t}\, d\sigma(t),\qquad z\in\mathbb{C}^{+}.\]
We will denote the above measure $\nu$ by $\nu_{\boxplus}^{\gamma,\sigma}$.
The following result is from \cite{BW in muli limit thm}. We reproduce
the short proof here because we actually require inequalities (2.1)
and (2.2). 

\begin{lem}
Consider a sequence $\{ r_{n}\}_{n=1}^{\infty}\subset\mathbb{R}$
and two triangular arrays $\{ z_{nk}:\, n\geq1,\,1\leq k\leq k_{n}\}$,
$\{ w_{nk}:\, n\geq1,\,1\leq k\leq k_{n}\}$ of complex numbers. Assume
that
\begin{enumerate}
\item $\Im w_{nk}\geq0$, for $n\geq1$ and $1\leq k\leq k_{n}$. 
\item \[
z_{nk}=w_{nk}(1+\varepsilon_{nk}),\]
where\[
\varepsilon_{n}=\max_{1\leq k\leq k_{n}}\left|\varepsilon_{nk}\right|\]
converges to zero as $n\rightarrow\infty$.
\item There exists a positive constant $M$ such that for sufficiently large
$n$,\[
\left|\Re w_{nk}\right|\leq M\Im w_{nk}.\]

\end{enumerate}
Then the sequence $\{ r_{n}+\sum_{k=1}^{k_{n}}z_{nk}\}_{n=1}^{\infty}$
converges if and only if the sequence $\{ r_{n}+\sum_{k=1}^{k_{n}}w_{nk}\}_{n=1}^{\infty}$
converges. Moreover, two sequences have the same limit.

\end{lem}
\begin{proof}
The assumptions on $\{ z_{nk}\}_{n,k}$ and $\{ w_{nk}\}_{n,k}$ imply
\begin{equation}
\left|\left(r_{n}+\sum_{k=1}^{k_{n}}z_{nk}\right)-\left(r_{n}+\sum_{k=1}^{k_{n}}w_{nk}\right)\right|\leq2(1+M)\varepsilon_{n}\left(\sum_{k=1}^{k_{n}}\Im w_{nk}\right),\label{eq:2.1}\end{equation}
and\begin{equation}
(1-\varepsilon_{n}-M\varepsilon_{n})\left(\sum_{k=1}^{k_{n}}\Im w_{nk}\right)\leq\left|\sum_{k=1}^{k_{n}}\Im z_{nk}\right|,\label{eq:2.2}\end{equation}
for sufficiently large $n$. If the sequence $\{ r_{n}+\sum_{k=1}^{k_{n}}z_{nk}\}_{n=1}^{\infty}$
converges to a complex number $z$, (2.2) implies that $\{\sum_{k=1}^{k_{n}}\Im w_{nk}\}_{n=1}^{\infty}$
is bounded, and then (2.1) shows that the sequence $\{ r_{n}+\sum_{k=1}^{k_{n}}w_{nk}\}_{n=1}^{\infty}$
also converges to $z$. Conversely, if $\{ r_{n}+\sum_{k=1}^{k_{n}}w_{nk}\}_{n=1}^{\infty}$
converges to $z$, then the sequence $\{\sum_{k=1}^{k_{n}}\Im w_{nk}\}_{n=1}^{\infty}$
is bounded and hence by (2.1) the sequence $\{ r_{n}+\sum_{k=1}^{k_{n}}w_{nk}\}_{n=1}^{\infty}$
converges to $z$ as well.
\end{proof}

\section{Proof of the Main Result}

Given an infinitesimal triangular array $\{\mu_{nk}:n\geq1,\,1\leq k\leq k_{n}\}\subset\mathcal{M}$,
define constants\[
a_{nk}=\int_{\left|t\right|<1}t\, d\mu_{nk}(t),\]
and measures $\overline{\mu}_{nk}$ by\[
d\overline{\mu}_{nk}(t)=d\mu_{nk}(t+a_{nk}).\]
Note that $\max_{1\leq k\leq k_{n}}\left|a_{nk}\right|\rightarrow0$
as $n\rightarrow\infty$, and this implies that $\{\overline{\mu}_{nk}\}_{n,k}$
is also an infinitesimal array. Define the analytic function\[
f_{nk}(z)=z^{2}\left[G_{\overline{\mu}_{nk}}(z)-\frac{1}{z}\right],\qquad z\in\mathbb{C}^{+},\]
and the real-valued function\[
b_{nk}(y)=\int_{\left|t\right|\geq1}a_{nk}\, d\mu_{nk}(t)+\int_{\left|t\right|\geq1}\frac{(t-a_{nk})y^{2}}{y^{2}+(t-a_{nk})^{2}}\, d\mu_{nk}(t),\qquad y\geq1.\]
Observe that $\Im f_{nk}(z)\leq0$ if $\Im z>0$, and $f_{nk}(z)=o(\left|z\right|)$
as $z\rightarrow\infty$ nontangentially.

The following lemma is related with a calculation in Section 4 of
\cite{CG2}.

\begin{lem}
For $y\geq1$, we have for all $n\in\mathbb{N}$, \[
\left|\Re\left[f_{nk}(iy)-b_{nk}(y)\right]\right|\leq2\left|\Im f_{nk}(iy)\right|,\]
and for sufficiently large $n$, \[
\left|\Re f_{nk}(iy)\right|\leq(3+6y)\left|\Im f_{nk}(iy)\right|,\]
where $1\leq k\leq k_{n}$.
\end{lem}
\begin{proof}
Note that for $y\geq1$, and $n\in\mathbb{N}$, \[
f_{nk}(iy)=\int_{-\infty}^{\infty}\frac{(t-a_{nk})y^{2}}{y^{2}+(t-a_{nk})^{2}}\, d\mu_{nk}(t)-i\int_{-\infty}^{\infty}\frac{(t-a_{nk})^{2}y}{(t-a_{nk})^{2}+y^{2}}\, d\mu_{nk}(t).\]
Moreover, since $\int_{\left|t\right|<1}(t-a_{nk})\, d\mu_{nk}(t)=\int_{\left|t\right|\geq1}a_{nk}\, d\mu_{nk}(t)$,
we have\begin{eqnarray*}
\left|\Re\left[f_{nk}(iy)-b_{nk}(y)\right]\right| & = & \left|\int_{\left|t\right|<1}\left[\frac{(t-a_{nk})y^{2}}{y^{2}+(t-a_{nk})^{2}}\, d\mu_{nk}(t)-(t-a_{nk})\right]\, d\mu_{nk}(t)\right|\\
 & = & \left|\int_{\left|t\right|<1}\frac{-(t-a_{nk})^{3}}{y^{2}+(t-a_{nk})^{2}}\, d\mu_{nk}(t)\right|\leq2\left|\Im f_{nk}(iy)\right|.\end{eqnarray*}
Note that the infinitesimality of the family $\{\mu_{nk}\}_{n,k}$,
implies that there exists $N\in\mathbb{N}$ such that $\left|a_{nk}\right|\leq1/2$,
for all $n\geq N$, $1\leq k\leq k_{n}$. Therefore, for $n\geq N$,
\begin{align*}
\left|\int_{\left|t\right|\geq1}a_{nk}\, d\mu_{nk}(t)\right| & \leq\left|a_{nk}\right|(1+4y^{2})\int_{\left|t\right|\geq1}\frac{(t-a_{nk})^{2}}{y^{2}+(t-a_{nk})^{2}}\, d\mu_{nk}(t)\\
 & \leq\frac{1+4y^{2}}{y}\left|\Im f_{nk}(iy)\right|\leq(1+4y)\left|\Im f_{nk}(iy)\right|,\end{align*}
and since $2x^{2}\geq\left|x\right|$ when $\left|x\right|\geq\frac{1}{2}$,
we have\begin{align*}
\left|\int_{\left|t\right|\geq1}\frac{(t-a_{nk})y^{2}}{y^{2}+(t-a_{nk})^{2}}\, d\mu_{nk}(t)\right| & \leq2y\left|\Im f_{nk}(iy)\right|.\end{align*}
Hence the second inequality follows. 
\end{proof}
\begin{lem}
Given $\beta\geq1$, let $\Gamma_{\alpha,\beta}$ be a truncated cone
where all the functions $\phi_{\mu_{nk}}$ are defined, and let $\{ c_{n}\}_{n=1}^{\infty}$
be a sequence of real numbers. 
\begin{enumerate}
\item For any $y>\beta$, the sequence $\{ c_{n}+\sum_{k=1}^{k_{n}}\phi_{\mu_{nk}}(iy)\}_{n=1}^{\infty}$
converges if and only the sequence $\{ c_{n}+\sum_{k=1}^{k_{n}}\left[a_{nk}+f_{nk}(iy)\right]\}_{n=1}^{\infty}$
converges. Moreover, the two sequences have the same limit.
\item If \[
L=\sup_{n\geq1}\sum_{k=1}^{k_{n}}\int_{-\infty}^{\infty}\frac{t^{2}}{1+t^{2}}\, d\overline{\mu}_{nk}(t)<\infty,\]
then $c_{n}+\sum_{k=1}^{k_{n}}\phi_{\mu_{nk}}(iy)=o(y)$ uniformly
in $n$ as $y\rightarrow\infty$ if and only if $c_{n}+\sum_{k=1}^{k_{n}}\left[a_{nk}+f_{nk}(iy)\right]=o(y)$
uniformly in $n$ as $y\rightarrow\infty$.
\end{enumerate}
\end{lem}
\begin{proof}
Fix $y>\beta$. From Proposition 2.3 applied to $\{\overline{\mu}_{nk}\}_{n,k}$,
we have\begin{equation}
-\phi_{\mu_{nk}}(iy)+a_{nk}=-\phi_{\overline{\mu}_{nk}}(iy)=-f_{nk}(iy)\cdot(1+v_{nk}(iy)),\label{eq:3.1}\end{equation}
where the functions \[
v_{n}(iy)=\max_{1\leq k\leq k_{n}}\left|v_{nk}(iy)\right|\]
converges to zero as $n\rightarrow\infty$. Then (1) follows from
Lemma 3.1 and Lemma 2.4 by choosing $z_{nk}=-\phi_{\mu_{nk}}(iy)+a_{nk}$,
$w_{nk}=-f_{nk}(iy)$ and $r_{n}=-c_{n}-\sum_{k=1}^{k_{n}}a_{nk}$.

We next prove (2). From (3.1), we have\[
z_{nk}^{\prime}(iy)=w_{nk}^{\prime}(iy)\cdot(1+v_{nk}(iy)),\]
where $z_{nk}^{\prime}(iy)=-\phi_{\mu_{nk}}(iy)+a_{nk}+b_{nk}(y)+b_{nk}(y)v_{nk}(iy)$,
and $w_{nk}^{\prime}(iy)=-f_{nk}(iy)+b_{nk}(y)$. Lemma 3.1 implies
that for all $n\in\mathbb{N}$, and $1\leq k\leq k_{n}$, \[
\left|\Re w_{nk}^{\prime}(iy)\right|\leq2\left|\Im f_{nk}(iy)\right|,\qquad y>\beta.\]
Since $\max_{1\leq k\leq k_{n}}\left|a_{nk}\right|\rightarrow0$ as
$n\rightarrow\infty$, there exists $N\in\mathbb{N}$ such that $\left|a_{nk}\right|\leq1/2$,
for all $n\geq N$, $1\leq k\leq k_{n}$. Therefore, for $n\geq N$,
and $y>\beta\geq1$,\begin{eqnarray*}
\sum_{k=1}^{k_{n}}\left|b_{nk}(y)\right| & \leq & \sum_{k=1}^{k_{n}}\int_{\left|t\right|\geq1}\frac{1}{2}\, d\mu_{nk}(t)+y\sum_{k=1}^{k_{n}}\int_{\left|t\right|\geq1}\left|\frac{(t-a_{nk})y}{y^{2}+(t-a_{nk})^{2}}\right|\, d\mu_{nk}(t)\\
 & \leq & (1+y)\sum_{k=1}^{k_{n}}\int_{\left|t\right|\geq1}\frac{1}{2}\, d\mu_{nk}(t)\leq5y\sum_{k=1}^{k_{n}}\int_{\left|t\right|\geq1}\frac{1}{5}\, d\mu_{nk}(t)\\
 & \leq & 5y\sum_{k=1}^{k_{n}}\int_{\left|t\right|\geq1}\frac{(t-a_{nk})^{2}}{1+(t-a_{nk})^{2}}\, d\mu_{nk}(t)\leq5yL.\end{eqnarray*}
Since $v_{n}(iy)=o(1)$ uniformly in $n$ as $y\rightarrow\infty$,
by decreasing the cone we may assume that $v_{n}(iy)<1/6$, for all
$y>\beta$, and $n\in\mathbb{N}$. Define $r_{n}^{\prime}(y)=-c_{n}-\sum_{k=1}^{k_{n}}a_{nk}-\sum_{k=1}^{k_{n}}b_{nk}(y)$.
Replacing $r_{n}$, $z_{nk}$, $w_{nk}$ by $r_{n}^{\prime}$, $z_{nk}^{\prime}$,
and $w_{nk}^{\prime}$ respectively in inequalities (2.1) and (2.2),
we deduce that\[
\left|\left(\sum_{k=1}^{k_{n}}\phi_{\mu_{nk}}(iy)\right)-\left(\sum_{k=1}^{k_{n}}\left[a_{nk}+f_{nk}(iy)\right]\right)\right|\leq\left|\sum_{k=1}^{k_{n}}\Im f_{nk}(iy)\right|+5yLv_{n}(iy),\]
 and\[
\frac{1}{2}\left|\sum_{k=1}^{k_{n}}\Im f_{nk}(iy)\right|\leq\left|\sum_{k=1}^{k_{n}}\Im\phi_{\mu_{nk}}(iy)\right|+5yLv_{n}(iy),\]
for all $n\geq N$, and $y>\beta$. Hence the result follows from
the facts that $v_{n}(iy)=o(1)$ uniformly in $n$ as $y\rightarrow\infty$,
and that (2) holds uniformly for $n$ in a finite subset of $\mathbb{N}$.
\end{proof}
Fix a real number $\gamma$ and a finite positive Borel measure $\sigma$
on $\mathbb{R}$. 

\begin{thm}
For an infinitesimal array $\{\mu_{nk}\}_{n,k}\subset\mathcal{M}$
and a sequence $\{ c_{n}\}_{n=1}^{\infty}\subset\mathbb{R}$, the
following statements are equivalent:
\begin{enumerate}
\item The sequence $\mu_{n1}\boxplus\mu_{n2}\boxplus\cdots\boxplus\mu_{nk_{n}}\boxplus\delta_{c_{n}}$
converges weakly to $\nu_{\boxplus}^{\gamma,\sigma}$;
\item The sequence $\mu_{n1}*\mu_{n2}*\cdots*\mu_{nk_{n}}*\delta_{c_{n}}$
converges weakly to $\nu_{*}^{\gamma,\sigma}$;
\item The sequence of measures\[
d\sigma_{n}(t)=\sum_{k=1}^{k_{n}}\frac{t^{2}}{1+t^{2}}\, d\overline{\mu}_{nk}(t)\]
 converges weakly on $\mathbb{R}$ to $\sigma$, and the sequence
of numbers \[
\gamma_{n}=c_{n}+\sum_{k=1}^{k_{n}}\left[a_{nk}+\int_{-\infty}^{\infty}\frac{t}{1+t^{2}}\, d\overline{\mu}_{nk}(t)\right]\]
converges to $\gamma$ as $n\rightarrow\infty$. 
\end{enumerate}
\end{thm}
\begin{proof}
The equivalence of (2) and (3) is classical (see \cite{GK,Petrov}).
We will prove the equivalence of (1) and (3). Assume that (1) holds.
By Theorem 2.2, there exist $\alpha>0$ and $\beta\geq1$ such that
$\phi_{\mu_{nk}}$ are defined on $\Gamma_{\alpha,\beta}$, and we
have\[
\lim_{n\rightarrow\infty}\phi_{\mu_{n1}\boxplus\mu_{n2}\boxplus\cdots\boxplus\mu_{nk_{n}}\boxplus\delta_{c_{n}}}(iy)=\phi_{\nu_{\boxplus}^{\gamma,\sigma}}(iy),\qquad y>\beta.\]
Since \[
\phi_{\mu_{n1}\boxplus\mu_{n2}\boxplus\cdots\boxplus\mu_{nk_{n}}\boxplus\delta_{c_{n}}}(z)=c_{n}+\sum_{k=1}^{k_{n}}\phi_{\mu_{nk}}(z),\qquad z\in\Gamma_{\alpha,\beta},\]
we have \[
\lim_{n\rightarrow\infty}\left(c_{n}+\sum_{k=1}^{k_{n}}\phi_{\mu_{nk}}(iy)\right)=\phi_{\nu_{\boxplus}^{\gamma,\sigma}}(iy),\]
for all $y>\beta$, and $c_{n}+\sum_{k=1}^{k_{n}}\phi_{\mu_{nk}}(iy)=o(y)$
uniformly in $n$ as $y\rightarrow\infty$. By Lemma 3.2, \begin{equation}
\lim_{n\rightarrow\infty}\left(c_{n}+\sum_{k=1}^{k_{n}}\left[a_{nk}+f_{nk}(iy)\right]\right)=\phi_{\nu_{\boxplus}^{\gamma,\sigma}}(iy),\qquad y>\beta.\label{eq:3.2}\end{equation}
Note that for $z\in\mathbb{C}^{+}$, $n\in\mathbb{N}$,\begin{eqnarray*}
z^{2}\left[G_{\overline{\mu}_{nk}}(z)-\frac{1}{z}\right] & = & \int_{-\infty}^{\infty}\frac{tz}{z-t}\, d\overline{\mu}_{nk}(t)\\
 & = & \int_{-\infty}^{\infty}\frac{t}{1+t^{2}}\, d\overline{\mu}_{nk}(t)+\int_{-\infty}^{\infty}\left[\frac{tz}{z-t}-\frac{t}{1+t^{2}}\right]\, d\overline{\mu}_{nk}(t)\\
 & = & \int_{-\infty}^{\infty}\frac{t}{1+t^{2}}\, d\overline{\mu}_{nk}(t)+\int_{-\infty}^{\infty}\left[\frac{1+tz}{z-t}\right]\frac{t^{2}}{1+t^{2}}\, d\overline{\mu}_{nk}(t).\end{eqnarray*}
We conclude that \begin{equation}
c_{n}+\sum_{k=1}^{k_{n}}\left[a_{nk}+f_{nk}(z)\right]=\gamma_{n}+\int_{-\infty}^{\infty}\frac{1+tz}{z-t}\, d\sigma_{n}(t).\label{eq:3.3}\end{equation}
Since $\Im f_{nk}(z)\leq0$ for $z\in\mathbb{C}^{+}$, $\{ c_{n}+\sum_{k=1}^{k_{n}}\left[a_{nk}+f_{nk}(z)\right]\}_{n=1}^{\infty}$
is a normal family of analytic functions in $\mathbb{C}^{+}$, and
from (3.2) the sequence has pointwise limit $\phi_{\nu_{\boxplus}^{\gamma,\sigma}}(z)$
for all $z=iy$, $y>\beta$. It is an easy application of the Montel
Theorem that (3.2) holds uniformly on compact subsets of $\mathbb{C}^{+}$.
Hence (3.3) and the formula of $\phi_{\nu_{\boxplus}^{\gamma,\sigma}}$
imply, at $z=i$, that \begin{eqnarray*}
\lim_{n\rightarrow\infty}\sigma_{n}(\mathbb{R}) & = & \lim_{n\rightarrow\infty}\int_{-\infty}^{\infty}\frac{1+t^{2}}{1+t^{2}}\, d\sigma_{n}(t)\\
 & = & \lim_{n\rightarrow\infty}-\Im\left(c_{n}+\sum_{k=1}^{k_{n}}\left[a_{nk}+f_{nk}(i)\right]\right)\\
 & = & \lim_{n\rightarrow\infty}-\Im\phi_{\nu_{\boxplus}^{\gamma,\sigma}}(i)\\
 & = & \sigma(\mathbb{R}).\end{eqnarray*}
Thus, \[
L=\sup_{n\geq1}\sigma_{n}(\mathbb{R})=\sup_{n\geq1}\sum_{k=1}^{k_{n}}\int_{-\infty}^{\infty}\frac{t^{2}}{1+t^{2}}\, d\overline{\mu}_{nk}(t)<\infty.\]
By Lemma 3.2, this implies that $c_{n}+\sum_{k=1}^{k_{n}}\left[a_{nk}+f_{nk}(iy)\right]=o(y)$
uniformly in $n$ as $y\rightarrow\infty$. For $y>\beta$, $n\in\mathbb{N}$,
note that\[
\frac{1}{2}\sigma_{n}(\{\left|t\right|\geq y\})\leq\int_{-\infty}^{\infty}\frac{1+t^{2}}{y^{2}+t^{2}}\, d\sigma_{n}(t)=-\frac{1}{y}\Im\left(c_{n}+\sum_{k=1}^{k_{n}}\left[a_{nk}+f_{nk}(iy)\right]\right).\]
Since $c_{n}+\sum_{k=1}^{k_{n}}\left[a_{nk}+f_{nk}(iy)\right]=o(y)$
uniformly in $n$ as $y\rightarrow\infty$, we conclude that $\{\sigma_{n}\}_{n=1}^{\infty}$
is a tight family. Let $\sigma^{\prime}$ be a weak cluster point
of $\{\sigma_{n}\}_{n=1}^{\infty}$ and a subsequence $\{\sigma_{n_{j}}\}_{j=1}^{\infty}$
converges weakly to $\sigma^{\prime}$. Hence, for any $z=x+iy\in\Gamma_{\alpha,\beta}$,
we have\begin{eqnarray*}
\int_{-\infty}^{\infty}\frac{y}{(x-t)^{2}+y^{2}}(1+t^{2})\, d\sigma^{\prime}(t) & = & -\lim_{j\rightarrow\infty}\Im\left(c_{n_{j}}+\sum_{k=1}^{k_{n_{j}}}\left[a_{n_{j}k}+f_{n_{j}k}(x+iy)\right]\right)\\
=-\Im\phi_{\nu_{\boxplus}^{\gamma,\sigma}}(x+iy) & = & \int_{-\infty}^{\infty}\frac{y}{(x-t)^{2}+y^{2}}(1+t^{2})\, d\sigma(t).\end{eqnarray*}
Therefore, the Poisson integrals of the measures $(1+t^{2})\, d\sigma^{\prime}(t)$
and $(1+t^{2})\, d\sigma(t)$ are identical since they coincide on
an open subset of $\mathbb{C}^{+}$. Thus, $\sigma^{\prime}=\sigma$.
Since the tight family $\{\sigma_{n}\}_{n=1}^{\infty}$ has a unique
weak cluster point, they must converge weakly to $\sigma$. Moreover,
we deduce that $\lim_{n\rightarrow\infty}\gamma_{n}=\gamma$ from
(3.2) and (3.3). 

To prove the converse, assume the sequence of measures $\{\sigma_{n}\}_{n=1}^{\infty}$
converges weakly to $\sigma$ and the sequence $\{\gamma_{n}\}_{n=1}^{\infty}$
converges to $\gamma$ as $n\rightarrow\infty$. Then, $\{\sigma_{n}\}_{n=1}^{\infty}$
is a tight family and in particular\[
L=\sup_{n}\sigma_{n}(\mathbb{R})<\infty.\]
From Lemma 2.1 and the infinitesimality of the array $\{\mu_{nk}\}_{n,k}$,
there exists a truncated cone $\Gamma_{\alpha^{\prime},\beta^{\prime}}$
with $\beta^{\prime}\geq1$ such that all $\phi_{\mu_{nk}}$ are defined
in $\Gamma_{\alpha^{\prime},\beta^{\prime}}$. Combine the inequality
\[
\left|\frac{1+ity}{iy-t}\right|\leq y,\qquad t\in\mathbb{R},\, y\geq1.\]
 with (3.3) to obtain\[
\lim_{n\rightarrow\infty}\left(c_{n}+\sum_{k=1}^{k_{n}}\left[a_{nk}+f_{nk}(iy)\right]\right)=\phi_{\nu_{\boxplus}^{\gamma,\sigma}}(iy),\qquad y>\beta^{\prime}.\]
 Hence by Lemma 3.2, \[
\lim_{n\rightarrow\infty}\left(c_{n}+\sum_{k=1}^{k_{n}}\phi_{\mu_{nk}}(iy)\right)=\phi_{\nu_{\boxplus}^{\gamma,\sigma}}(iy),\qquad y>\beta^{\prime}.\]
Also, note that for any $M>0$ and $y>\beta^{\prime}\geq1$, we have\begin{eqnarray*}
\frac{1}{y}\left|c_{n}+\sum_{k=1}^{k_{n}}\left[a_{nk}+f_{nk}(iy)\right]\right| & \leq & \frac{\left|\gamma_{n}\right|}{y}+\frac{1}{y}\int_{-\infty}^{\infty}\left|\frac{1+ity}{iy-t}\right|\, d\sigma_{n}(t)\\
 & \leq & \frac{\left|\gamma_{n}\right|}{y}+\frac{1}{y}\int_{\left|t\right|<M}\frac{1+\left|t\right|y}{\sqrt{y^{2}+t^{2}}}\, d\sigma_{n}(t)+\sigma_{n}(\{\left|t\right|\geq M\})\\
 & \leq & \frac{\left|\gamma_{n}\right|}{y}+\frac{L(1+My)}{y^{2}}+\sigma_{n}(\{\left|t\right|\geq M\}).\end{eqnarray*}
 Therefore, it follows from the convergence of $\{\gamma_{n}\}_{n=1}^{\infty}$
and the tightness of the family $\{\sigma_{n}\}_{n=1}^{\infty}$ that
$c_{n}+\sum_{k=1}^{k_{n}}\left[a_{nk}+f_{nk}(iy)\right]=o(y)$ uniformly
in $n$ as $y\rightarrow\infty$ . By Lemma 3.2 again, $c_{n}+\sum_{k=1}^{k_{n}}\phi_{\mu_{nk}}(iy)=o(y)$
uniformly in $n$ as $y\rightarrow\infty$. Statement (1) now follows
from Theorem 2.2.
\end{proof}

\end{document}